\documentclass{conm-p-l}
\usepackage{amssymb,amsmath,amsthm,amscd } %,float}, amsaddr} %,color,mathtools, pdflscape}

\usepackage[T1]{fontenc}

\usepackage{booktabs, topcapt}

\usepackage{hyperref}

\hypersetup{
  colorlinks   = true, %Colours links instead of ugly boxes
  urlcolor     = blue, %Colour for external hyperlinks
  linkcolor    = blue, %Colour of internal links
  citecolor   = red %Colour of citations
}

\usepackage[msc-links]{amsrefs}

\usepackage{mathtools}

\def\Z{\mathbb Z}
\def\Q{\mathbb Q}
\def\R{\mathbb R}

\def\P{\mathbb P}
\def\wP{\mathbb{WP}}

\newcommand\M{\mathcal M}
\def\O{\mathcal O}
\newcommand\p{\mathfrak p}

\newcommand\x{\mathbf{x}}
\newcommand\y{\mathbf{y}}
\newcommand\w{\mathbf{w}}

\def\<{\langle}
\def\>{\rangle}

\newcommand\A{\mathbb A}
\newcommand{\Orb}{\mathit{Orb }}

\def\X{\mathcal X}

\def\l{\lambda}
\def\w{\mathfrak w}
\def\x{\mathbf x}

\def\i{\mathfrak i}

\DeclareMathOperator\wgcd{\mathit{wgcd }}   %   weighted gcd
\DeclareMathOperator\wh{\mathfrak{h}}   %  weighted height

\DeclareMathOperator\awh{\mathfrak{\tilde h} }   % absolute    weighted height

\newtheorem{thm}{Theorem}
\newtheorem{prop}{Proposition}
\newtheorem{lem}{Lemma}
\newtheorem{cor}{Corollary}
\newtheorem{rem}{Remark}
\newtheorem{exa}{Example}
\newtheorem{prob}{Problem}

\theoremstyle{definition}
\newtheorem{defi}[thm]{Definition}

\theoremstyle{remark}

\def\bK{\bar K}
\def\bQ{\overline \Q}

\def\awgcd{\overline{\wgcd }}

\def\q{\mathfrak q}

\begin{document}
\title{Computing heights on weighted projective spaces}

\author{Jorgo Mandili}
\address{Research Institute of Sciences and Technology,   Vlorë, Albania.}
\email{mandili@risat.org}
\author{Tony Shaska}
\address{Department of Mathematics and Statistics,  Oakland University,    Rochester, MI 48309. }
\email{shaska@oakland.edu}

\dedicatory{Dedicated to  Mehmet Likaj, \\ on the occasion of his 70th birthday}

\begin{abstract}
In this  note we extend the  concept  height on projective spaces to that of weighted height on weighted projective spaces and show how such a height can be computed.
We prove some of the basic properties of the weighted height and show how it can be used to study hyperelliptic curves over $\mathbb Q$. 
Some examples are provided from the weighted moduli space of binary sextics and  octavics. 
\end{abstract}


\maketitle

\def\J{\mathfrak J}

%\tableofcontents
%********************************************
\section{Introduction}

Let  $\wP^n_w (K)$ be the weighted projective space of weights $\w=(w_0, \dots , w_n)$ over a field $K$ of characteristic zero.  Is there a way to measure the "size" of points in  $\wP^n_w (K)$ similar to the height on the projective space $\P^n (K)$?  The answer comes from \cite{b-sh} where the concept of the height was defined for weighted projective spaces. How does one compute the height of a point in $\wP^n_w (K)$?  Moreover, how do we get a tuple $\x= (x_0, \dots x_n)$ such that it is the minimal representative for the point $\p =[x_0, \cdots , x_n ] \in \wP^n (K)$?   In this short paper we explore these questions.

The motivation for considering  the above question comes from the theory of hyperelliptic   or superelliptic curves.  The isomorphism classes of a genus $g\geq 2$ hyperelliptic curve $C: \ y^2 z^{2g}=f(x, z)$ correspond to the tuple of generators of the ring of invariants $S (2, 2g+2)$  of binary forms evaluated at the binary form $f(x, z)$.  Such ring of invariants    is a weighted projective space.  Hence, determining a canonical minimal  tuple for any point in  $\wP^n_w (K)$ would give a one to one correspondence between the isomorphism classes of curves and such minimal tuples. We illustrate briefly with the genus 2 curves. 

In \cite{data} we created a database of isomorphism classes of genus 2 curves defined over $\Q$.  Every such isomorphism class was identified uniquely by a set of absolute invariants $(i_1, i_2, i_3)$; see \cite{data} for details.  These invariants are defined in terms of the Igusa invariants $J_2, J_4, J_6, J_{10}$. Why not  identify the curve with the tuple $(J_2, J_4, J_6, J_{10})$ instead of  $(i_1, i_2, i_3)$?  If we do so then we have determine how to pick the smaller size tuple for any point $\p= [ J_2, J_4, J_6, J_{10}]$ and how to do this in a canonical way. The goal of this paper is to address such issues for any weighted projective space.

In this paper we define  a  \textit{normalization} of points   $\p \in \wP^n_w (\Q)$ which is the representing tuple of $\p$ with smallest coefficients.  We show that this normalization is unique up to multiplication by a primitive $d$-th root of unity, where $d= \gcd (w_0, \dots , w_n)$ and is unique when $\wP^n_w (\Q)$ is well-formed.  The   height of a point  $\p \in \wP^n_{\w} (\Q)$ is the weighted absolute value of coordinates of $\p$, when $\p$ is normalized. 
%We generalized such concept of height over any number field $K$.  

We also define the \textit{absolutely normalized} tuples which is a normalization over the algebraic closure  $\bar \Q$.  This is a normalization by multiplying by scalars which are allowed to be in $\bar \Q$.   The height of an absolutely normalized tuple is called an \textit{absolute height} for analogy with the terminology in \cite{height}. In other words the \textit{absolute height} of a point $\p \in \wP^n_{\w} (\Q)$  is the  the weighted absolute value of coordinates of $\p$, when $\p$ is absolutely normalized. 
%We generalized the concept of absolute height over any number field $K$.  

The paper is organized as follows.  In Section~\ref{sect-2} we give a brief introduction to weighted projective spaces $\wP^n_w (K)$.  A standard reference here is \cite{dol}. We consider both \textit{well-formed} and not \textit{well-formed} weighted projective spaces.  For any point $\x = (x_0, \dots , x_n)  \in (\Z)^{n+1} \setminus \{ 0 \} $,  we define the \textit{weighted greatest common divisor} $\wgcd (\x)$ as  the product of all primes $p \in \Z$ such that for all $i=0, \dots , n$, we have $p^i \mid x_i$.  We will call a point $\p \in \wP_{\w}^n (\Q)$ \textbf{normalized} if it has $\wgcd (\p) =1$. 
%We also generalize the normalization of points in $\wP_{\w}^n (K)$ over the algebraic closure $\bar k$.  

In Section~\ref{sect-3} we follow \cite{b-sh} and  define the \textit{weighted projective height} on $\wP^n_w (\Q)$ and show that this is well-defined.  We prove a version of the Northcott's theorem for the weighted projective height and determine for what conditions on integers $w_o, \dots w_n$ the normalized tuple is unique.    
Analogously we extend the definitions and results over the algebraic closure.   We show how to determine all twists of a given point $\p \in \wP_{\w}^n (\Q)$ of height $h \leq \wh (\p)$.  When the set of weights is $\w=(1, \dots , 1)$, then the weighted projective space is simply the projective space $\P^n (K)$ and our weighted moduli height becomes the usual height on $\P^n (K)$ as defined in \cite{bombieri-book}.   
 
The notion of weighted height and absolute weighted height is used in  \cite{scott}, to study the weighted moduli space of binary sextics and in \cite{w-3} to study the weighted moduli space of binary octavics. Both cases lead to creating databases of genus 2 or genus 3 hyperelliptic curves with small absolute moduli height.  For connections of weighted projective spaces and the algebraic curves  or other topics on databases of hyperelliptic curves the reader can check \cite{data}. We give some examples for genus 2 curves and genus 3 hyperelliptic curves, which were the main motivation behind this paper. It remains to be seen if there are any explicit relations between the weighted moduli height,  moduli height, and height as in \cite{height}.

The concept of weighted height in weighted projective spaces, surprisingly seems unexplored before. The only reference we could find was the unpublished report in \cite{deng} which defines the function \textit{Size} similarly to our height with different motivations. Our goal in writing this short note was to simply provide a brief introduction to heights in weighted projective spaces. 
We assume the reader is familiar with the concept of height in projective spaces as in \cite{bombieri-book} and \cite{silv-book}. 
A more detailed study of heights on weighted projective spaces over any number field  is given in \cite{b-sh, F-S}. \\

%\smallskip

\noindent \textbf{Notation}  The algebraic closure of a field $K$ is denoted by  $\bK$.   
For  an algebraically number field $K$ we denote by  $\O_K$ its ring of integers and by $M_K$   the set of all absolute values in $K$.  
A point $\p \in \wP_{\w}^n (K)$ is denoted by $\p = \left[ x_0 : x_1 : \dots : x_n  \right]$ and the tuple of coordinates $\x= (x_0, x_1, \dots , x_n)$. By a "curve" we always mean the isomorphism class of a smooth, irreducible curve.  
%*************************
\section{Weighted greatest common divisors}

Let $q_0$, \dots , $q_n$ be positive integers.  A set of weights is called the ordered tuple 
\[
\w=(q_0, \dots , q_n).
\] 
Denote by  $r=\gcd (q_0, \dots , q_n)$ the greatest common divisor of $q_0, \dots , q_n$.  A \textit{weighted integer tuple} is a tuple $\x = (x_0, \dots x_n ) \in \Z^{n+1}$ such that to each coordinate $x_i$ is assigned the weight $q_i$.  We multiply weighted tuples by scalars $\lambda \in \Q$ via 
\[ \lambda \star (x_0, \dots , x_n) = \left( \l^{q_0} x_0, \dots , \l^{q_n} x_n   \right) \]

For an ordered   tuple of integers  $\x=(x_0, \dots , x_n) \in \Z^{n+1}$, whose coordinates are not all zero, the \textbf{weighted greatest common divisor with respect to the set of weights} $\w$ is the largest integer $d$ such that 
\[ d^{q_i} \, \mid \, x_i, \; \; \text{for all }    i=0, \dots n.\]
 We will denote the weighted greatest common divisor by $\wgcd (x_0, \dots x_n)$.   A tuple $\x=(x_0, \dots , x_n)$ with $\wgcd (\x)=1$ is called \textbf{normalized}. 

The \textbf{absolute weighted greatest common divisor} of  a tuple $\x=(x_0, \dots , x_n)$ with respect to the set of weights $\w=(q_0, \dots , q_n)$ is the largest integer $d$ such that 
\[ d^{\frac {q_i} r} \, \mid \, x_i, \; \; \text{for all }    i=0, \dots n.\]
 We will denote the weighted greatest common divisor by $\awgcd (x_0, \dots x_n)$.   A tuple $\x$ with $\awgcd (\x)=1$ is called \textbf{absolutely normalized}.
 
\begin{exa}
Consider the set of weights $\w=(2, 4, 6, 10)$ and a tuple 
\[ \x = \left(   3 \cdot 5^2 , 3^2 \cdot 5^4, 3^3 \cdot 5^6, 3^5 \cdot 5^{10} \right) \in \Z^4.\]
Then,  $\wgcd (\x) = 5$ and \; $\awgcd (\x) = 3 \cdot 5$.
\end{exa} 

We summarize in the following lemma.

\begin{lem}
For any weighted integral tuple $\x= (x_0, \dots , x_n) \in \Z^{n+1}$, the tuple
\[   \y = \frac 1 {\wgcd (\x) } \star \x, \]
is integral and normalized. Moreover, the tuple 
\[    \bar \y = \frac 1 {\awgcd (\x) } \star \x, \]
is also integral and absolutely normalized. If $\gcd (q_0, \dots , q_n)=1$, then  $\wgcd (\x) = \awgcd (\x)$. 
\end{lem}

%\proof
The proof is a direct consequence of the definition.
%\qed

Let us attempt to extend these definitions to any ring of integers.  Let $K$ be any number field and $\O_K$ its ring of integers. Consider the set of weights $\w$ as above and an integer tuple $\x \in \O_K^{n+1}$. For any $\alpha \in \O_K$, the ideal generated by alpha is denoted by $(\alpha)$.

The \textbf{weighted greatest common divisor ideal} is defined as 
\[ \J (\x) = \sum_{(x_i) \subset \left( \p^{q_i} \right)}  (\p) \]
over all primes $\p$ in  $\O_K$.  The \textbf{absolute weighted greatest common divisor ideal} is defined as 
\[ \J (\x) = \sum_{(x_i) \subset \left( \p^{\frac {q_i} r} \right)}  (\p) \]
over all primes $\p$ in  $\O_K$.   In general the weighted greatest common divisor is defined for all Dedekind domains or more generally for all GCD-domains.

%********************
%\newpage
\section{Weight projective spaces}\label{sect-2}
%*******************************************
Let $K$ be a field of characteristic zero and  $(q_0, \dots , q_n) \in \Z^{n+1}$ a fixed tuple of positive integers called \textbf{weights}.  
 Consider the action of $K^\star = K \setminus \{0\}$ on $\A^{n+1} (K)$ as follows
\[ \lambda \star (x_0, \dots , x_n) = \left( \l^{q_0} x_0, \dots , \l^{q_n} x_n   \right) \]
for $\l\in K^\ast$. 
The quotient of this action is called a \textbf{weighted projective space} and denoted by   $\wP^n_{(q_0, \dots , q_n)} (K)$. 
%The space $\P(1, \dots , 1)$ is the usual projective space. 
It is the projective variety $Proj \left( K [x_0,...,x_n] \right)$ associated to the graded ring $K [x_0, \dots ,x_n]$ where the variable $x_i$ has degree $q_i$ for $i=0, \dots , n$. 

We denote greatest common divisor of $q_0, \dots , q_n$ by $\gcd (q_0, \dots , q_n)$.    The space $\wP_w^n$ is called \textbf{well-formed} if   
\[ \gcd (q_0, \dots , \hat q_i, \dots , q_n)   = 1, \quad \text{for each } \;  i=0, \dots , n. \]
While most of the papers on weighted projective spaces are on well-formed spaces, we do not assume that here.   We will denote a point $\p \in \wP_w^n (K)$ by $\p = [ x_0 : x_1 : \dots : x_n]$. 

Let $K$ be a number field and $\O_K$ its ring of integers. The group action $K^\star$ on $\A^{n+1} (K)$ induces a group action of $\O_K$ on $\A^{n+1} (K)$.         
By  $\Orb (\p)$ we denote the $\O_K$-orbit in $\A^{n+1} (\O_K)$ which contains $\p$. 
For any point $\p =[x_0 : \dots : x_n] \in \wP_w^n (K)$   we can assume, without loss of generality, that $\p =[x_0 : \dots : x_n] \in \wP_w^n (\O_K)$. The  height for weighted projective spaces will be defined in the next section.  

%Next we define some auxiliary concepts which are very helpful in computations. 
For the rest of this section we assume $K=\Q$.   
For the tuple $\x = (x_0, \dots , x_n) \in \Z^{n+1}$ we define   the \textbf{weighted greatest common divisor} with respect to the absolute value $| \, \cdot \, |_v$, denoted by $\wgcd_v (\x)$,   
\[ \wgcd_v (\x) :=  \prod_{  \stackrel {d^{q_i} \mid x_i} {d\in \Z}    } | d |_v \]
as the product of all divisors  $d\in \Z$ such that for all $i=0, \dots , n$, we have $d^i \mid x_i$.  
We will call a point $\p \in \wP_{\w}^n (\Q)$ \textbf{normalized} if $\wgcd (\p) =1$. 
%The  \textbf{weighted greatest common divisor} of a tuple $\x =(x_0 , \dots , x_n) \in \A^{n+1} (\Q)$ is defined to be
%
%\[ \wgcd (\x) :=  \prod_{  \stackrel {p^{q_i} \mid x_i} {p\in \Z}    } | \, p \, | \]
%
%for all primes $p$ in $\Z$.  

\begin{defi}
We will call a point $\p \in \wP_{\w}^n (\Q)$ a \textbf{normalized point} if the weighted greatest common divisor of its coordinates is 1.  
\end{defi}

\begin{lem}\label{lem1} 
Let $\w=(q_0, \dots , q_n)$ be a set of weights and $d=\gcd (q_0, \dots , q_n)$. 
For any point $\p \in \wP_{\w}^n (\Q)$, the point
\[ \q = \frac 1 {\wgcd (\p)} \star \p \]
is normalized.  Moreover, this normalization  is unique up to a multiplication by a $d$-root of unity.
\end{lem}

\proof  Let $\p = [x_0 : \dots , x_n ]  \in \wP_w^n (\Q)$ and $\p_1 = [\alpha_0 : \dots : \alpha_n ]$ and $\p_2 = [ \beta_0 : \dots : \beta_n]$ two different normalizations of $\p$. Then exists non-zero $\l_1, \l_2 \in \Q$   such that 
\[ 
\p = \l_1 \star \p_1 = \l_2 \star  \p_2,
\]
or in other words 
\[ 
(x_0, \dots , x_n) = \left( \l_1^{q_0} \alpha_0 , \dots , \l_1^{q_i} \alpha_i , \dots \right) =  \left( \l_2^{q_0} \beta_0 , \dots , \l_2^{q_i} \beta_i , \dots \right).
\]
Thus,
\[ 
  \left(   \alpha_0 , \dots ,   \alpha_i , \dots , \alpha_n \right) =  \left( r^{q_0} \beta_0 , \dots , r^{q_i} \beta_i , \dots , r^{q_n} \beta_n \right).
\]
for $r= \frac {\l_2} {\l_1} \in K$. Thus, $r^{q_i}=1$ for all $i=0, \dots , n$.  Therefore,  $r^d=1$.  This completes the proof.   
\qed

Thus we have the following:

\begin{cor}
For any point $\p=[x_0 : \cdots : \x_n] \in \wP_{\w}^n (\Q)$, if the greatest common divisors of non-zero coordinates is 1,  
%$\gcd (q_0, \dots , q_n)=1$, 
then the normalization of   $\p$ is unique. 
\end{cor}

Here is an example which illustrates the Lemma.

\begin{exa}
Let $\p = [x_0, x_1, x_2, x_3]  \in \wP_{(2,4,6,10)}^3 (\Q)$ be a normalized point.  
Hence, 
\[ \wgcd (x_0, x_1, x_2, x_3 ) =1.\]
Since $d = \gcd (2, 4, 6, 10)=2$, then we can take $r$ such that   $r^2=1$.  Hence, $r=\pm 1$.   Therefore,  the point
\[ (-1)  \star \p = [-x_0 : x_1 : -x_2 : -x_3 ] \]
is also be normalized.  

However, if $\p = [x_0, x_1, x_2, x_3]  \in \wP_{(1, 2, 3, 5)}^3 (\Q)$ is normalized then it is unique, unless some of the coordinates are zero. For example the points $[0, 1, 0, 0]$ and $[0, -1, 0, 0]$ are equivalent and both normalized. 
\end{exa}

Next we give two examples, which were the main motivation behind this note. 

%*****************
\begin{exa}[Weighted projective space of binary sextics]
The ring of invariants of binary sextics is generated by the \textbf{basic arithmetic invariants}, or as they sometimes called,  \textbf{Igusa invariants} $(J_2 , J_4 , J_6 , J_{10})$ as defined in \cite{Ig-60}. 
Two genus 2 curves $\X$ and $\X^\prime$ are isomorphic if and only if there exists $\l \in K^\ast$ such that 
\[ J_{2i} (\X) = \l^{2i} \, J_{2i} (\X^\prime), \quad \text{for } \quad i=1, 2, 3, 5. \]
We take the set of weights $\w=(2, 4, 6, 10)$ and considered the weighted projective space $\wP_{(2, 4, 6, 10)} (\Q)$.  Thus, the invariants of a sextic define a point in a weighted projective space $[J_2 : J_4 : J_6 : J_{10}] \in \wP_{\w} (\Q)$ and every genus 2 curve correspond to a point in $\wP_{\w}^3 (\Q) \setminus \{ J_{10}\neq 0 \}$.  There is a bijection between
\[ \phi : \; \; \wP_{(2, 4, 6, 10)}^3 \setminus \{ J_{10} \neq 0 \} \to \M_2, \]
with $\phi$ provided explicitly in \cite[Theorem 1]{univ-gen-2}. 
\end{exa}
Using the notion of a normalized point as above we have the following:

\begin{cor} Normalized points in $\wP_{(2, 4, 6, 10)}^3 (\Q)$ occur in pairs.  In other words, 
for every normalized point $\p= \left[J_2, J_4, J_6, J_{10}\right]$, there is another normalized point   
$\p^\prime = \left[-J_2, J_4, -J_6, -J_{10}\right]$ equivalent to $\p$.  Moreover, $\p$ and $\p^\prime$ are isomorphic over the Gaussian integers. 
\end{cor}

\proof
Let $\X$ be a genus 2 curve with equation $y^2=f(x)$ and $[J_2, J_4, J_6, J_{10}]$ its corresponding invariants.  The transformation $x \mapsto \sqrt{-1} \cdot x$ with give a curve $\X^\prime$ with invariants $[- J_2 : J_4 : - J_6 : -J_{10}]$ and the same weighted moduli height. 

If two weighted moduli points have the same minimal absolute height, then they differ up to a multiplication by a unit. Hence, 
\[ [J_2^\prime : J_4^\prime : J_6^\prime : J_{10}^\prime] = [d^2\cdot J_2 : d^4\cdot J_4 : d^6\cdot J_6 : d^{10}\cdot J_{10}] \]
such that $d^2$ is a unit.  Then, $d^2=\pm 1$. Hence, $d=\sqrt{-1}$. 
\qed

So unfortunately for any genus 2 curve we have two corresponding normalized points $[ \pm J_2, J_4, \pm J_6, \pm J_{10}]$. In \cite{scott} this problem is solve by taking always the point $[ | J_2|, J_4, \pm J_6, \pm J_{10}]$ or by considering the space $\wP_{(1, 2, 3, 5)}^3 (\Q)$ instead.

%**************************************************************************************************
\begin{exa}[Weighted projective space of binary octavics]

Every irreducible, smooth, hyperelliptic genus 3 curve has  equation $y^2 z^6=f(x, z)$, where $f(x, z)$ is a binary octavic with non-zero discriminant.  The ring of invariants of binary octavics is generated by invariants $J_2, \dots , J_8$, which satisfy an algebraic equation as in \cite[Thm.~6]{gen-3}. Two genus 3 hyperelliptic curves $\X$ and $\X^\prime$ 
%in Weierstrass form, given by equations  
%\[  C:   Z^2=f(X, Y)   \textit{    and   } C^\prime:  z^2=g(X, Y) \]
%
are isomorphic over a field  $K$ if and only if there exists some $\l \in k\setminus \{ 0\}$ such that 
\[ J_i (\X) = \l^i J_i(\X^\prime), \textit{   for   }  \,\,  i=2, \dots , 7.   \]  
%
%and $J_2, \dots J_8$ satisfy the Eq.~\eqref{shaska}.
There is another invariant $J_{14}$ given in terms of $J_2, \dots J_7$ which is the discriminant of the binary octavic.

Hence, there is a bijection between the hyperelliptic locus in the moduli space of genus 3 curves and the weighted projective space $\wP_{(2,3,4,5,6,7)}^5 (K) \setminus \{ J_{14} \neq 0\}$. Since $d = \gcd (2,3,4,5,6,7)=1$ then we have:
\end{exa}

\begin{cor}
For every genus 3 hyperelliptic curve $\X$, defined over a field $K$, the corresponding normalized point 
\[ \p = [J_2 : J_3 : J_4 : J_5 : J_6 : J_7] \in \wP_{(2,3,4,5,6,7)}^5 (K)\]
is unique. 
\end{cor}

\begin{exa}\label{gen-3-ex-1}
Consider   the curve $y^2=x^8-1$.  The moduli point in $\wP_{\w}^5 (\Q)$ is 
\[ \p = \left[-    2^{3} \cdot 5 \cdot 7, 0,   2^{10} \cdot 7^{4}, 0,  2 ^{15} \cdot 7 ^{6}, 0,-  2 ^{19} \cdot 5 \cdot 7^{8} \right]   \]
Then,    $\wgcd (\x) =   \frac 1  2$. 
Hence, the point $\p$  normalized becomes 
\[ \frac 1 2 \star \p = \left[-    2  \cdot 5 \cdot 7, 0,   2^6 \cdot 7^{4}, 0,  2^9 \cdot 7^{6}, 0, -  2^{11} \cdot 5 \cdot 7^{8} \right].  \]
\end{exa}

In \cite{w-3} we use such normalized points to create a database of genus 3 hyperelliptic curves defined over $\Q$. 

%*********************************
\subsection{Absolutely normalized points}
For  any point $\p = [x_0 : \dots : x_n ] \in \wP_{\w}^n (\Q)$ we may assume that $x_i \in \Z$ for $i=0, \dots , n$
and define 
\[ 
\overline{\wgcd} (\p) = \prod_{\l \in \bar \Q,  \,  \l^{q_i} | x_i} | \l |
\]
as the product of all $\l \in \bar \Q$, such that for all $i=0, \dots , n$, $\l^i \in \Z$ and $\l^i | x_i$.  
A point $\p = [x_0 : \dots : x_n ] \in \wP_{\w}^n (\Q)$ is called \textbf{absolutely normalized} or \textbf{normalized over $\bar \Q$} if $\awgcd (\p)=1$. 

%*****************************************
%\subsection{Over number fields}

%We generalize over any number field $K$ as follows.   For  any point $\p = [x_0 : \dots : x_n ] \in \wP_{\w}^n (K)$ we define 
%
%\[ \overline{\wgcd} (\p) = \prod_{\l \in \bK,  \l^{q_i} | x_i} \l \]
%
% Equivalently for  any point $\p = [x_0 : \dots : x_n ] \in \wP_{\w}^n (K)$ we define  the following fractional ideal  associated to $\p$ 
%
% \[ \a^{-1} (\p)= \{ \lambda \in \bK \; : \; \lambda * (x_0, \dots , x_n) \in \O_K \} \]
%
%Then $\a^{-1} (\p)= \< \overline{\mathfrak \delta} \>$.  We have $\overline{\wgcd} (\p)= \overline{\mathfrak \delta}$. 

\begin{defi}
A point $\p = [x_0 : \dots : x_n ] \in \wP_{\w}^n (\Q)$ is called \textbf{absolutely normalized} or \textbf{normalized over the algebraic closure} if $\awgcd (\p)=1$. 
\end{defi}

\begin{lem}\label{lem2}
For any point $\p = [x_0 : \dots : x_n ] \in \wP_{\w}^n (\Q)$ its normalization over the algebraic closure
\[ {\bar \p} = \frac 1 {\awgcd (\p)} \star \p \]
is unique up to a multiplication by a $d$-th root of unity.  
\end{lem}

\proof  Let $\p = [x_0 : \dots , x_n ]  \in \wP_w^n (\Q)$ and $\p_1 = [\alpha_0 : \dots : \alpha_n ]$ and $\p_2 = [ \beta_0 : \dots : \beta_n]$ two different normalizations of $\p$ over $\bQ$. Then exists non-zero $\l_1, \l_2 \in \bQ$   such that 
\[ 
\p = \l_1 \star \p_1 = \l_2 \star  \p_2,
\]
or in other words 
\[ 
(x_0, \dots , x_n) = \left( \l_1^{q_0} \alpha_0 , \dots , \l_1^{q_i} \alpha_i , \dots \right) =  \left( \l_2^{q_0} \beta_0 , \dots , \l_2^{q_i} \beta_i , \dots \right).
\]
Thus,
\[ 
  \left(   \alpha_0 , \dots ,   \alpha_i , \dots , \alpha_n \right) =  \left( r^{q_0} \beta_0 , \dots , r^{q_i} \beta_i , \dots , r^{q_n} \beta_n \right).
\]
for $r= \frac {\l_2} {\l_1} \in \bQ$. Thus, $r^{q_i}=1$ for all $i=0, \dots , n$.  Therefore,  $r^d=1$.  This completes the proof.   
\qed

Two points  $\p$ and $\q$  in $\wP_{\w}^n (\Q)$  are called \textbf{twists } of each other if they are equivalent  in $\wP_{\w}^n (\bQ)$  but  $\Orb_{\Q} (\p )$ is not the same as  $\Orb_{\Q} ( \q )$.  Hence, we have the following. 
\begin{lem}
Let $\p$ and $\p^\prime$ be normalized points in $\wP_{\w}^n (\Q)$.  Then  
 $\p$ and $\p^\prime$ are twists of each other if and only if there exists $\l \in \bQ^\star$ such that $\l \star \p = \p^\prime$.
\end{lem}

Next we see another example from genus 2 curves.

\begin{exa}
Let $\X$ be the genus two curve with equation $y^2=x^6-1$ and $J_2, J_4, J_6$, and $J_{10}$ its Igusa   invariants. Then the isomorphism class of $\X$ is determined by the point
$ \p = [240, 1620, 119880, 46656] \in \wP_{(2, 4, 6, 10)}^3 (\Q)$.    Thus, % its weighted moduli point is
\[ \p = [240, 1620, 119880, 46656] = [ 2^4\cdot 3 \cdot 5   ; \,   2^2 \cdot 3^4 \cdot 5 ; \, 2^3 \cdot 3^4 \cdot 5  \cdot 37 ; \, 2^6\cdot 3^6 ].\]
Therefore,
\[ 
\begin{split}
& \wgcd (240, 1620, 119880, 46656) = 1 \\
&  \awgcd (240, 1620, 119880, 46656) = \sqrt{6}.
\end{split}
\]
Hence, $\p$ is normalized  but not absolutely normalized.     
The point $\p$  has   twists, %for values of $d= 2, 3, 6$, which are 
\[
\begin{split}
\p_1 & = \frac 1 {\sqrt{2}} \star \p_1 =[120, 405, 14985, 1458] =  [ 2^3\cdot 3 \cdot 5 \; : \;    3^4 \cdot 5 \; :\;   3^4 \cdot 5 \cdot  37\; :\; 2\cdot 3^6 ],\\
\p_2 & = \frac 1 {\sqrt{3}} \star \p_1 = [80, 180, 4440, 192]= [ 2^4 \cdot 5 \; : \;  2^2 \cdot 3^2 \cdot 5 \; :\; 2^3 \cdot 3 \cdot 5 \cdot  37\; :\; 2^6\cdot 3 ],\\
\end{split}
\]
and the absolutely normalized point of $\p$ which is 
\[ {\bar \p} = \frac 1 {\sqrt{6}} \star \p_1 = [40, 45, 555, 6] = [ 2^3\cdot 5, \, 3^2 \cdot 5, \, 3\cdot 5 \cdot 37 , \, 2\cdot 3]\]
Notice that ${\bar \p}$ has only one twist 
\[  {\bar \p}^\prime = [ -2^3\cdot 5, \, 3^2 \cdot 5, \, -3\cdot 5 \cdot 37 , \, -2\cdot 3]\]
which is also normalized. 
\end{exa} 

We can do better even with the genus 3 curve from Example~\ref{gen-3-ex-1}.

\begin{exa}
The normalized  moduli point in $\wP_{\w}^5 (\Q)$ the curve $y^2=x^8-1$ is  
\[ \frac 1 2 \star \p = \left[-    2  \cdot 5 \cdot 7, 0,   2^6 \cdot 7^{4}, 0,  2^9 \cdot 7^{6}, 0, -  2^{11} \cdot 5 \cdot 7^{8} \right].  \]
Then,
$ \awgcd (\p) = \frac {\i} {\sqrt{14} }$,  for $ \i^2 =-1$.
Then its absolutely normalized form is
\[ {\bar \p} = \left[  5 ,  0,   2^4 \cdot 7^2, 0,  2^6 \cdot 7^3, 0, -  2^7 \cdot 5 \cdot 7^4 \right].  \]
\qed
\end{exa}

In the next section we will introduce some measure of the magnitude of points in weighted moduli spaces $\wP_{\w}^n (K)$ and show that the process of normalization and absolute normalization lead us to the representation of points in $\wP_{\w}^n (K)$ with smallest possible coordinates. 

%***************************
%\newpage
\section{Heights on the weighted projective spaces}\label{sect-3}

%\subsection{Reduction }

%\begin{defi} 
%We say that a point $\p \in \wP_w^n (\O_K)$ is \textbf{normalized} in $\O_K$ if there is no non-unit  $d \in O_K$ such that $d^i \mid x_i$, for all $i=0, \dots , n$. 
%\end{defi}

%In this section we define a \textit{height} or \textit{magnitude} on the weighted projective spaces.  

Let $K$ be an algebraic number field and $[K : \Q]=n$ and its ring of integers $\O_K$.  With $M_K$ we denote the set of all absolute values in $K$. For $v \in M_K$, the \textbf{local degree at $v$}, denoted $n_v$ is  $n_v =[K_v:\Q_v]$,   where $K_v, \Q_v$ are the completions with respect to $v$. As above 
 $\wP^n (K)$ is the projective space with weights $w=(q_0, \dots, q_n)$, and  $\p \in \wP^n(K)$ a point with   coordinates $\p=[x_0, \dots , x_n]$, for  $x_i \in K$. The \textbf{multiplicative heigh}t of $\p$ is defined as follows
\[\wh_K(\p) := \prod_{v \in M_K} \max\left\{\frac{}{}|x_0|_v^{n_v/q_0} , \dots, |x_n|_v^{n_v/q_n}\right \}\]
Let $\p=[x_0, \dots , x_n] \in \wP^n(\Q)$ with weights $w = (q_0, \dots, q_n)$.  It is clear that $\p$ will have a representative $[y_0, \dots, y_n]$ such that $y_i \in \Z$ for all $i$ and $\wgcd(y_0, \dots, y_n)=1$.  With such representative for the coordinates of $\p$, the non-Archimedean absolute values give no contribution to the height, and we obtain
\[\wh_\Q (\p)=  \max_{0 \leq j \leq n}\left\{\frac{}{}|x_j|^{1/q_j}_\infty \right\}\]
% 

\iffalse 

We define the \textbf{height} or  \textbf{magnitude} of a  point  $\p = [x_0 : \dots : x_n ] \in \wP_w^n (\O_K)$ as  
%
\begin{equation}\label{height} 
\wh (\p)=  \frac 1 {\wgcd(\p)} \,  \prod_{v \in M_K} \max  \left\{ |x_0|_v^{1/q_0},  \dots , |x_n|_v^{1/q_n}  \right\},   
\end{equation}
%
where $M_K$ is the set of all norms in $K$.  Let $\p_0 = [y_0 : \dots : y_n ] $ be the normalization of $\p$.  Then obviously 
%
\[ \wh (\p) = \wh (\p_0) = \prod_{v \in M_K} \max  \left\{ |y_0|_v^{1/q_0},  \dots , |y_n|_v^{1/q_n}  \right\}.\]
%

\fi

So for a tuple  $\x = (x_0 : \dots : x_n )$ the height of the corresponding point $\p = [\x]$  is 
\[ 
\wh (\p) = \frac 1 {\wgcd (\x) } \max  \left\{ |x_0|^{1/q_0},  \dots , |x_n|^{1/q_n}  \right\}. 
\]
%Notice that $\wh_{\infty}$ is as follows
%
%\begin{equation}\label{h-inf}
% \wh_{\infty} (P) = \frac 1 {\displaystyle \prod_{ p^{q_i} |  x_i }  p  }     \, \max \{  |x_0|^{1/q_0},  \dots , |x_n|^{1/q_n}|  \}  
% \end{equation}
%
%for all primes $p \in \O_K$. 

We combine some of the properties of $\wh (\p)$ in the following: 
\begin{prop}\label{prop-1}  
Then the  following are true: 

i)  The   function  $\wh  : \wP_{\w}^n (\Q) \to \R$  is well-defined. 

ii) A normalized point $\p = [x_0 : \dots : x_n]  \in \wP_{\w}^n (\Q)$ is the point with  smallest coordinates in its orbit $\Orb (\p)$. 

iii) For any constant $c > 0$ there are only finitely many points $\p \in \wP_w^n (\Q)$ such that    $\wh (\p) \leq c$. 
\end{prop}

\proof
i)   It is enough to show that two normalizations of the same point $\p \in \wP_{\w}^n (\Q)$ have the same height. Let $\p$ and $\q$ be such normalizations.  Then from Lemma~\ref{lem1} we have $\p = r \star \q$, where $r^d=1$. Thus,
\[ \wh (\p) = \wh ( r \star \q ) = |r| \cdot \wh (\q) = \wh (\q). \]

%Let $\p=(x_0, \dots , x_n)$ and  $\l \in \O_K$ such that 
%
%\[ \l = p_1^{n_1} \cdots p_r^{n_r}. \]
%
%Then,  $\l*\p=\left( \l^{q_0} x_0, \dots , \l^{q_n} x_n  \right)$ and  we have 
%
%\[  \begin{split}
% \wh ( \l*\p ) & = \frac 1 { \mathfrak d }    \, \, \max_{0\leq i \leq n} \{ | \l^{q_i} \cdot x_i |^{\frac 1 {q_i}} \}    = \frac 1 {\l} \cdot \l \cdot \max_{0\leq i \leq n} \{ |   x_i |^{\frac 1 {q_i}} \} = \wh(\p), 
% \end{split} \]
%
%where $\mathfrak d = \displaystyle \prod_{p^{q_i} | \l^{q_i} x_i } p $, for all $p\in \O_K$.

ii) This is obvious from the definition. 

iii) Let $\p \in \wP_{\w}^n (\Q)$.   It is enough to count only normalized points    $\p= [x_0 : \dots : x_n]  \in \wP_w^n (\Z)$ such that $\wh (\p) \leq c$.  For every coordinate $x_i$ there are only finitely values in $\Z$ such that $|x_i|_v^{1/q_i} | \leq c$.  Hence, the result holds. 

\qed

Part iii) of the above  is the analogue of the Northcott's theorem in projective spaces.  
%We prove it only over $\Q$; see \cite{height} for a general statement over a number field $K$. 

\begin{rem}
If the set of weights $\w=(1, \dots 1)$ then $\wP_{\w}^n (\Q)$ is simply the projective space $\P^n (\Q)$ and the height $\wh (\p)$ correspond to the height of a projective point as defined in \cite{height}. 
\end{rem}

Let's see an example how to compute the height of a point.

\begin{exa}
Let $\p = (2^2, 2\cdot 3^4, 2^6 \cdot 3 , 2^{10} \cdot 5^{10}) \in \wP_{(2,4,6,10)}^3 (\Q)$. Notice that $\p$ is normalized, which implies that 
\[   
\wh (\p) =   \max \, \left\{ 2, 2^{1/4} \cdot 3,   3^{1/6} , 2  \cdot 5    \right\}        = 10
\]
However, the point $\mathfrak q = (2^2, 2^4 \cdot 3^4, 2^6 \cdot 3 , 2^{10} \cdot 5^{10}) \in \wP_{(2,4,6,10)}^3 (\Q)$ can be normalized to $(1, 3^4, 3, 5^{10})$ which has height
\[ \wh (\mathfrak q) =  \max \, \left\{ 1,   3,   3^{1/6} ,  5    \right\}    =  5. \]
%\qed
\end{exa}

A proof for the following will be provided in \cite{b-sh}.

\begin{lem}\label{lem_1}
Let $\p \in \wP^n (K)$ with weights $w = (q_0, \dots, q_n)$ and $L/K$ be a finite extension. Then,
\[\wh_L (P)= \wh_K (P)^{[L:K]}.\]
\end{lem}

%*************************************************

% can be written as 
%
%\begin{equation}
% \mathfrak d (\p) = \displaystyle \prod_{p^{q_i} \in \left( x_i \right)} p 
%\end{equation}
%
%where the product is taken over all primes $p \in \O_K$ such that $p^{q_i} \in \left( x_i \right)$, for $i=0, \dots , n$. 
%We define the \textbf{magnitude} of $\p = (x_0, \dots , x_n ) \in \wP_w^n (\O_K)$ as  
%
%\begin{equation}\label{height} 
%\wh (\p)= \frac 1 {| \mathfrak d |} \,  \prod_{v \in M_K} \max  \left\{ |x_0|^{1/q_0},  \dots , |x_n|^{1/q_n}  \right\}    
%\end{equation}

%********************************************
%\newpage
\subsection{Absolute heights}\label{sect-4}

   We  can define the height on $\wP^n(\overline \Q)$. The height of a point on $\wP^n(\overline \Q)$ is called the \textbf{weighted absolute (multiplicative)  height} and is the function 
\[
\begin{split}
\awh: \wP^n(\bar \Q) & \to [1, \infty)\\
\awh (\p)&= \wh_K (P)^{1/[K:\Q]},
\end{split}
\]
where  $\p \in \wP^n(K)$, for any $K$.   Then, the absolute weight height is given by 
\begin{equation}\label{abs-height} 
\awh_{\Q} (\p)= \frac 1 { \overline{\wgcd} (\p)     } \,   \max  \left\{ |x_0|^{1/q_0},  \dots , |x_n|^{1/q_n}  \right\}    
\end{equation}
%
%where $\mathfrak N ( \a (P) )$ is the nilradical of $\a (P)$. 
%\[ \h (\p)=  \frac  { \max \{ |a_1|^{1/q_1},  \dots , |a_n|^{1/q_n}| \} }  {\sum |a_i^{1/q_i}| }   .\]
%
 
Let's see an example which compares the height of a point with the absolute height. 

\begin{exa} Let $\p = [ 0: 2 : 0 : 0] \in \wP_{(2,4,6,10)}^3 (\Q)$.  Then $\p$ is normalized and therefore $\wh (\p) = 2$.  However, it absolute normalization is 
$ \q = \frac 1 {2^{1/4} } \star \p = [ 0 : 0 : 1 : 0]$.
Hence, $\awh (\p ) = 1$. 
\end{exa}

\begin{rem}
As a consequence of the above results it is possible to "sort" the points in $\wP_{\w}^n (\bK)$ according to the absolute height and even determine all the twists for each point when the weighted projective space is not well-formed.  This is used in   \cite{scott} to create a database of genus 2 curves and similarly in \cite{w-3} for genus 3 hyperelliptic curves. 
\end{rem}

The  \textbf{weighted absolute height}   of $\p=[\x]\in \wP_{\w}^n (K)$, where $\x=( x_0 : \dots : x_n )$, for any number field $K$, is     
\begin{equation}\label{abs-height} 
\awh_K (\p)= \frac 1 { \overline{\wgcd} (\x)     } \,  \prod_{v \in M_K} \max  \left\{ |x_0|^{1/q_0},  \dots , |x_n|^{1/q_n}  \right\}    
\end{equation}
%
%where $\mathfrak N ( \a (P) )$ is the nilradical of $\a (P)$. 
%\[ \h (\p)=  \frac  { \max \{ |a_1|^{1/q_1},  \dots , |a_n|^{1/q_n}| \} }  {\sum |a_i^{1/q_i}| }   .\]
%

The concept of weighted absolute height correspond to that of absolute height in \cite{height}.  In \cite{height} a curve with minimum absolute height has an equation with the smallest possible coefficients.  In this paper, the absolute height says that there is a representative tuple of $\p \in \wP_{\w}^n (K)$ with smallest magnitude of coordinates.

Then we have the following:

\begin{prop}\label{prop-2}  
Let $K$ be a number field and  $\O_K$ its ring of integers.  Then the  following are true: 

i)  The absolute height function  $\awh_K \, : \wP_{\w}^n (K) \to \R$  is well-defined. 

ii) $\awh (\p)$ is the minimum of  heights of all twists of $\p$. 

iii) For any constant $c>0$ there are only finitely many points $\p \in \wP_w^n (K)$ such that    $\awh (\p) \leq c$. 
\end{prop}

\proof Part ii) and iii) are obvious.  We prove part i). We have to show that two different normalizations over the algebraic closure have the same absolute height.  Let $\p$ and $\q$ be such normalizations.  Then from Lemma~\ref{lem2} we have $\p = r \star \q$, where $r^d=1$. Thus,
\[ \awh(\p) = \awh ( r \star \q ) = |r| \cdot \awh (\q) = \wh (\q). \]
%
%\[ \wh (\p) = \wh ( r \star \q ) = |r| \cdot \wh (\q) = \wh (\q). \]

This completes the proof.
\qed

For more details we direct the reader to \cite{b-sh}.
Let's revisit again our example from genus 2 curves. 

\begin{exa}
Let $\X$ be the genus two curve with equation $y^2=x^6-1$ and moduli  point
$ \p = [240, 1620, 119880, 46656] \in \wP_{(2, 4, 6, 10)}^3 (\Q)$. We showed that  $\p$ is normalized  and therefore has  height $\wh (\p) = 4 \sqrt{15}$.
Its absolute normalization is
\[ {\bar \p} = [40, 45, 555, 6] = [ 2^3\cdot 5, \, 3^2 \cdot 5, \, 3\cdot 5 \cdot 37 , \, 2\cdot 3]\]
Hence, the    absolute height is  $\awh (\p) =  2 \sqrt{10}$.
\end{exa} 

%****************************************
\section{Computing the weighted height}

Given a  point $\p \in \wP^n (K)$, how easy is it  to compute its weighted height $\wh_K (\p)$? From the previous section this would be equivalent to computing the weighted greatest common divisor $\wgcd (\x)$ for a point $\x \in \O_K^{n+1}$, such that $\p=[\x]$. 
There are issues to be resolved when computing over $\O_K$, so for the purposes of this paper we continue to assume $K=\Q$.  

Computing $\wgcd (\x)$ is equivalent with factoring every coordinate over $\Z$.  Hence, this approach is not very effective for points with large coordinates. Hence, the main part of concern for any algorithm of computing the weighted height of a point $\p \in \wP^n (\Q)$ is the \textit{normalization} of a point in $\wP^n (\Q)$.  We have implemented this algorithm in SageMath and it has been used in \cite{scott} and \cite{w-3} to create databases of binary sextics and binary octavics of small weighted height.  It works well for small heights $\wh$. Recall that the size of the coordinates for any $\p \in \wP^3_{2, 4, 6, 10} (\Q)$  and $\wh (\p) \leq c$ is $\leq c^{10}$. In general, for a point $\p \in \wP^n$ with maximal weight among coordinates $w$, the worst bound is $c^w$.

One of the main problems of arithmetic related to such heights is the following.  Consider $f(x, y)\in K[x, y]$ a binary form of degree $d\geq 2$.  In classical mathematics determining  conditions on the coefficients of $f(x, y)$ such that $f(x, y)$ has minimal discriminant has been well studied, but only understood for small $n$ (i.e., $n=2, 3$). However, the discriminant is only one of the invariants of the degree $d$ binary forms, so one of the coordinates of the corresponding point $\p$  in the weighted projective space $\wP_w (K)$.  The complete problem would be restated as:

\begin{prob}
Determine conditions on the coefficients of the binary form $f(x, y)$ such that the corresponding point $\p$  in the weighted projective space $\wP_w (K)$ is absolutely normalized. 
\end{prob}

This seems out of reach for any degree $n>3$. However, our algorithm suggested above would work well over $\Q$ in finding a binary form $g$ equivalent to $f$ over $\Q$ with such that the corresponding point in $\wP^n (\Q)$ is normalized.  For more on this problem see \cite{b-sh-2}. 

In \cite{b-sh-2} we suggest an algorithm to compute equations of hyperelliptic or superelliptic curves which correspond to a normalized point in the weighted projective space.  This algorithm is an extension of Tate's algorithm of elliptic curves \cite{ta-75} and methods suggested in \cite{rachel}. 

There is another problem that comes from the analogy with the discriminant.  A classical problem as determining the number of curves with bounded discriminant (or a good bound for such number of curves), becomes now the problem of determining a good bound for the number of curves with bounded weighted moduli point.  It is unclear of any such good bounds, since now we don't want to estimate the number of tuples with bounded weighted height, but the number of equivalence classes of such tuples. Number that must be significantly less than the number of tuples. Some heuristically data for the space $\wP^3_{2, 4, 6, 10} (\Q)$ of binary sextics is given in \cite{scott}.

%*****************************

%\section{Concluding remarks}

Heights on weighted projective spaces, surprisingly have not been explored before.  There is an unpublished preprint by A. Deng (1998)    with the intention of counting the rational points in weighted projective spaces; see \cite{deng}. 
%The author in \cite{deng} defines a function called \textit{Size} which is somewhat simi
This is the first article where the concept of the height is defined in weighted projective spaces. A full account of heights in weighted projective spaces and their properties is intended in \cite{b-sh}.

%*****************
%\newpage
\nocite{*}
\bibliographystyle{amsplain} 

\bibliography{ref}{}

\end{document}